\documentclass[preprint,12pt]{elsarticle}
\usepackage[left=3cm,right=3cm]{geometry}
\usepackage[latin1]{inputenc}
\usepackage{amsmath}
\usepackage{amsfonts}
\usepackage{amssymb}
\usepackage{makeidx}
\usepackage{graphicx}
\usepackage{amsthm}
\usepackage{mathrsfs}
\numberwithin{equation}{section}
\usepackage{ragged2e}
\usepackage{tikz-cd}
\usepackage{csquotes}
\usepackage{natbib}
\newtheorem{thm}{Theorem}

\newtheorem{definition}[thm]{\quad Definition}

\usepackage{graphics}
\usepackage[all,cmtip]{xy}
\usepackage{graphicx}
\usepackage[toctitles]{titlesec}
\usepackage{amsfonts}
\usepackage{latexsym}
\usepackage{amssymb}
\usepackage{amsmath}
\usepackage{mathtools}
\usepackage{color}
\newtheorem{theorem}{Theorem}[section]
\newtheorem{lemma}[theorem]{Lemma}

\newtheorem{proposition}[theorem]{Proposition}
\newtheorem{corollary}[theorem]{Corollary}
\newtheorem{remark}[theorem]{Remark}

\journal{...}
\begin{document}
	
	\begin{frontmatter}
		
		\title{Automorphisms of Multiplicative Lie algebra Extensions}
		\author{Dev Karan Singh$^1$,  Shiv Datt Kumar$^2$ \vspace{0cm}}
		\address{$^{1,2}$Department of Mathematics\\ Motilal Nehru National Institute of Technology Allahabad, Prayagraj (UP) - $211004$, India, $^1$devkaransingh1811@gmail.com, dev.2020rma01@mnnit.ac.in, $^2$sdt@mnnit.ac.in}
		\begin{abstract}
			In this paper, we discuss the inducibility problem for automorphisms of multiplicative Lie algebra extensions and show that obstruction to the inducibility of pairs lies in the second cohomology group of multiplicative Lie algebras. We also establish the Wells type exact sequence for multiplicative Lie algebras, which relates automorphism groups with the second cohomology group of multiplicative Lie algebras.  
		\end{abstract}
	\end{frontmatter}
	\section{Introduction}
	The concept of multiplicative Lie algebra was introduced by G. J. Ellis \cite{GJ} for giving a non-abelian version of the Magnus and Witt theorem, which essentially states that five well known universal commutator identities generate all the commutator identities of weight $n$.
	A multiplicative Lie algebra is a group (possibly non-abelian) with an extra binary operation called multiplicative Lie product, which satisfies five identities similar to five well known universal commutator identities.
	The theory of multiplicative Lie algebra is interesting and distinct on its own. In the past few years authors have shown interest to develop the theory for multiplicative Lie algebra. In (\cite{MRS},\cite{FPW}) theory of nilpotency and solvability for multiplicative Lie algebras have been studied.
	Bak et al. in \cite{AGNM} constructed two homology theories for multiplicative Lie algebras and compared them with the homology theories of groups and Lie rings. They also introduced the Steinberg multiplicative Lie algebra of a unital ring and prove that it's the direct product of the Steinbereg group and Steinberg Lie ring. Lal and Upadhyay in \cite{RLS} discussed the extension theory for multiplicative Lie algebra and drived the Schur-Hopf formula for multiplicative Lie algebra in terms of second cohomology. They also discussed some applications to $K$- theory and cyclic homology. Donadze et al. in \cite{GNMA} introduced the concept of non-abelian tensor product and obtained a six term exact sequence for multiplicative Lie algebra and also proved a new version of Stalling's theorem.
	Pandey and Upadhyay in \cite{MSc} discussed the classification of multiplicative Lie products for finite groups. Pandey and Upadhyay in \cite{MSS} consider a short exact sequence $E(H,K) \equiv  \xymatrix{ e \ar[r] & H\ar[r]^{i} & G \ar[r]^{\beta} & K\ar[r] & e}$, where $H$ is an abelian group with trivial multiplicative Lie product, $G$ and $K$ are multiplicative Lie algebra. A short exact sequence  $E(H,K)$ is called a center extension if $H$ is contained in the center of $G$. They also discussed the extension theory for multiplicative Lie algebras and defined the second cohomology group for multiplicative Lie algebras.
	
	C. Wells in \cite{Wells} formulated an exact sequence for the automorphism group of a given group extension, providing a solution to the inducibility problem for a pair of automorphisms. In (\cite{Jin},\cite{Manoj}) authors further explored this sequence, focusing on studying the lifting and extensions of automorphisms in abelian extensions and established Wells-type exact sequences for specific cases. In \cite{VB_Msingh} Bardakov and Singh discussed a similar problem in the context of Lie algebras and developed Wells-type sequence for Lie algebras. In this paper, our primary objective is to address the inducibility problem for multiplicative Lie algebras, as discussed in Section 2, analogous to groups \cite{Wells} and Lie algebras \cite{VB_Msingh}. We also see that the obstruction for a pair of automorphisms to be inducible lies in the second cohomology group of the multiplicative Lie algebra.
	
	The rest of the paper is organized as follows. In Section 2, we provide a brief description of the second cohomology group for multiplicative Lie algebras, as outlined in \cite{MSS}. Subsequently, we introduce the notion of inducible and compatible pairs for multiplicative Lie algebras, establishing a necessary and sufficient condition for a pair $(\alpha, \eta)\in Aut(H)\times Aut(K)$ to be inducible in a center extension. In Section 3, we discuss the problem under what conditions an automorphism of an ideal of multiplicative Lie algebra can be extended to the multiplicative Lie algebra itself which fixes the ideal. Following the spirit of C. Wells \cite{Wells}, we also establish the Wells type exact sequence for multiplicative Lie algebras. As an application of these sequences, we discuss the inducibility of a pair of automorphisms. Additionally, we prove that for a split center extension, every compatible pair is inducible. We also obtain a necessary condition for the inducibility of the pair $(\alpha, \eta)\in Aut(H)\times Aut(K)$ when $K$ is a perfect multiplicative Lie algebra. In Section 4, we prove  that if the center extension $E(H,K)$ of multiplicative Lie algebras splits, then the exact sequences obtained in Section 3 also split.
	\section{Exact Sequences and Inducibility of a pair of Automorphisms}
	First, we give the definition of multiplicative Lie algebra and a brief introduction of the extension theory for multiplicative Lie algebras.
	\begin{definition}\label{multiplicative}
		A multiplicative Lie algebra is a triple $(G, \cdot, \star )$, where $G$ is group with respect to the binary operation $\cdot$ and satisfy the following conditions with respect to the binary operation $\star$
		\begin{enumerate}
			\item $x \star x = e$,
			\item $x\star y\cdot z = (x \star y)\cdot~^y(x \star z)$,
			\item  $x\cdot y \star z = ~^x(y \star z)\cdot (x \star z)$,
			\item $((x\star y ) \star ~^yz)\cdot((y\star z) \star ~^zx)\cdot((z \star x) \star ~^xy) = e$,
			\item $^z(x \star y) = ~^zx \star~^zy$, for all $x,y,z \in G$.
		\end{enumerate} 
		The binary operation $\star $ is called a multiplicative Lie product and $^zx$ denotes $zxz^{-1}$. 
	\end{definition}
	\noindent \textbf{Note:} A multiplicative Lie algebra is a generalization of groups and Lie algebra and every multiplicative Lie algebra need not be a Lie algebra but every Lie algebra can be seen as multiplicative Lie algebra by defining $\star$ to as Lie bracket of Lie algebra.
	
	Let $E(H,K)\equiv \xymatrix{e\ar[r] & H\ar[r]^{i} & G \ar[r]^{\beta} & K\ar[r] & e}$ be a center extension and $t: K \to G$ be a transversal of $E(H,K)$. Here $H$ is an ideal of $G$, i.e. $x \star h \in H$ for all $x \in G$ and $h \in H$.  We can also see that every element $g \in G$ can be uniquely written as $g = t(x)h$ for some $h \in H$ and $x \in K$. Also, for $x,y \in K$ and $h \in H$, we have $\beta(t(xy))=xy = \beta(t(x)t(y)), ~ \beta(t(x\star y)) = x\star y = \beta(t(x) \star t(y))$ and $\beta(t(x)\star h) = x \star e= e$. Thus we get maps $f: K \times K \to H, ~ h: K \times K \to H$ and $\Gamma_x : H \to H$ defind by $f(x,y) = t(x) t(y) t(xy)^{-1}, ~ h(x,y)= (t(x) \star t(y))(t(x\star y))^{-1}$ and $\Gamma_x(h) = t(x) \star h$ (It can be easily seen that $\Gamma_x \in End(H)$). Note that $End(H)$ is a multiplicative Lie algebra with operations defined as $(f \cdot g )(h )= f(h)\cdot g(h)$ and $(f \star g)(h) = f(g(h)) \cdot g(f(h^{-1}))$ for all $f,g \in End(H)$ and $h \in H$. Using the definition of multiplicative Lie algebra, we can establish the following properties of the maps $f,h$ and $\Gamma_x$.
	\begin{eqnarray}  
		\label{eq: 1}
		f(e,x) = f(x,e) =e\\
		\label{eq: 2}
		h(x,x)=h(x,e)=h(e,x)=e\\
		\label{eq: 3}
		f(x, y)f(xy, z) = f(y, z)f(x, yz)
	\end{eqnarray} \vspace{-1.2 cm}
	\begin{multline} 
		\label{eq: 4}
		f(y^{-1}, y(x \star z))^{-1}f(y(x \star z), y^{-1})f(x \star y, y(x \star z))h(x, y)h(x, z)
		= \Gamma_x(f(y, z))h(x, yz)
	\end{multline} \vspace{-1 cm}
	\begin{multline}
		\label{eq: 5}
		f(x^{-1}, x(y \star z))^{-1}f(x(y \star z), x^{-1})f(x(y \star  z), x \star z)h(x, z)h(y, z) = \Gamma_z(f(x, y)^{-1})h(xy, z)
	\end{multline} \vspace{-1 cm}
	\begin{multline}
		\label{eq: 6}
		\Gamma_x(f(z^{-1}, zy)^{-1}f(zy, z^{-1}))\Gamma_y(f(z^{-1}, zx)f(zx, z^{-1})^{-1})h(^zx,~ ^zy) \\
		= f(z^{-1}, z(x \star y))^{-1}f(z(x \star y), z^{-1})h(x, y)
	\end{multline} \vspace{-1.2 cm}
	\begin{multline}
		\label{eq: 7}
		\Gamma_{x\star y}(f(y^{-1}, yz)^{-1}f(yz, y^{-1}))\Gamma_{y\star z}(f(z^{-1}, zx)^{-1}
		· f(zx, z^{-1}))\Gamma_{z\star x}(f(x^{-1}, xy)^{-1}f(xy, x^{-1}))\\ \Gamma_z(h(x, y)^{-1})
		· \Gamma_x(h(y, z)^{-1})\Gamma_y(h(z, x)^{-1})h((x \star y),~^yz)h((y \star z),~^zx)h((z \star x),~^xy) \\f((x\star y)\star ~^yz,(y\star z)\star ~^zx)f(((z\star x)\star ~^xy)^{-1},((z\star x)\star ~^xy)))=1
	\end{multline}
	for all $x,y$ and $z \in K$. A triple $(f,h,\Gamma)$ is called a multiplicative center 2-cocycle of $K$ with coefficient in $H$, if $(f,h,\Gamma)$ satisfy the above properties and the set of all such triples is denoted by $Z^2_{ML(\Gamma)}(K,H)$ and its an abelian group with respect to the operation $(f,h,\Gamma)\cdot (f',h',\Gamma)=(ff',hh',\Gamma)$. A triple $(f,h, \Gamma)$ is called a multiplicative center 2-coboundary of $K$ with coefficient in $H$ if there exists a identity preserving map $\lambda:K \to H$ such that
	\begin{enumerate}
		\item $f(x,y) = \lambda(x) \lambda(y) \lambda(xy)^{-1}$
		\item $h(x,y) =  \Gamma_x(\lambda(y)) \Gamma_y(\lambda(x)^{-1}) \lambda(x \star y)^{-1}$ for all $x,y \in K$.
	\end{enumerate}
	The set of all such triples is denoted by $B^2_{ML(\Gamma)}(K,H)$ and its subgroup of $Z^2_{ML(\Gamma)}(K,H)$. The quotient group $H^2_{ML(\Gamma)}(K,H) = Z^2_{ML(\Gamma)}(K,H)/B^2_{ML(\Gamma)}(K,H)$ is called the second center cohomology of $K$ with coefficient in $H$. Two multiplicative centers $2$-cocycle are said to be cohomologous if they differ by a multiplicative center $2$-coboundary. 
	Let $\nu$ be a map from $K$ to $H$ and $\Gamma: K \to End(H)$ be a multiplicative Lie algebra homomorphism. Then $\nu$ is said to be mutiplicative center $1$-cocycle, if the following holds 
	\begin{enumerate}
		\item $\nu(xy) =\nu(x) \nu(y)$,
		\item $\nu(x \star y)= \Gamma_x(\nu(y)) \Gamma_y(\nu(x)^{-1}) $	for all $x,y \in K$.
	\end{enumerate}
	We denote the set of all such maps by $Z_{ML(\Gamma)}^1(K,H)$, which is an abelian group under the usual operation.\\
	
	\noindent \textbf{Note:} For more details of the extension theory of multiplicative Lie algebras, readers can see \cite{MSS}.	In the entire paper, we consider the center extension $E(H,K)$, in which, the ideal $H\subseteq Z(G)$ and multiplicative Lie product $\star$ on $H$ is trivial. Thus we have 
	\begin{equation}
		\label{equation_first}
		hg=gh, h\star k=e ~\text{and}~ h\star g \in H 
	\end{equation}
	for all $h,k\in H$ and $g\in G$.  We  will use throughout the paper properties mentioned in  \ref{equation_first} and Definition \ref{multiplicative} (as well as its properties) without explicitly mentioning them, especially  to  show that certain maps are homomorphisms. \\
	
	Now we will discuss the inducibility problem for the center extension of multiplicative Lie algebras.
	
	Let $E(H,K)\equiv \xymatrix{e\ar[r] & H\ar[r]^{i} & G \ar[r]^{\beta} & K\ar[r] & e}$ be a center extension and $t: K \to G$ be a transversal of $E(H,K)$. Define $Aut_H(G)  = \{ \phi \in Aut(G) : \phi(H) = H\}$. Then any $\phi \in Aut_H(G)$ induces a pair $(\alpha, \eta) \in Aut(H)\times Aut(K)$, which is given by $\alpha(h)=\phi(h)$ for all $h\in H$ and $\eta(x)=\beta(\phi(t(x)))$ for all $x \in K$.  
	Thus we can define a map \[\Pi : Aut_H(G) \to Aut(H) \times Aut(K), ~\text{given by}~ \Pi(\phi) = (\alpha, \eta).\]
	The homomorphism property of $\phi$ is easily apparent. However in the next lemma we see that $\eta$ is independent of the choice of transversal $t$.
	
	\begin{lemma}\label{zero_lemma}
		If $\phi \in \text{Aut}_H(G)$, then the map $\eta$ does not depend on the choice of the transversal $t$.
		
		\begin{proof}
			Let $t$ and $t'$ be two transversals of $E(H,K)$. Then $\beta(t(x))=x=\beta(t'(x))$ for all $x\in K$. This implies $t(x)t'(x)^{-1}\in \text{ker}(\beta)=H$. Thus $\phi(t(x)t'(x)^{-1})\in H = \text{ker}(\beta)$ for $\phi \in \text{Aut}_H(G)$. Therefore $\beta(\phi(t(x))\phi(t'(x)^{-1}))=e$ implies $\beta(\phi(t(x)))=\beta(\phi(t'(x)))$.
		\end{proof}
	\end{lemma}
	From the map $\Pi$, it follows that every $\phi \in \text{Aut}_H(G)$ induces a pair $(\alpha, \eta) \in \text{Aut}(H) \times \text{Aut}(K)$. The map $\Pi$ motivates us to introduce the following definition and raises the following question.
	\begin{definition}
		Let $\Pi : Aut_H(G) \to Aut(H) \times Aut(K)$ be a homomorphism. Then a pair $(\alpha, \eta) \in \text{Aut}(H) \times \text{Aut}(K)$ is called  inducible if there exists $\phi \in \text{Aut}_H(G)$ such that $\Pi(\phi) = (\alpha, \eta)$.
	\end{definition}
	\textbf{Question:} Under what conditions,  a pair of multiplicative Lie algebra automorphisms in $\text{Aut}(H) \times \text{Aut}(K)$ is inducible?
	
	In the next result, we provide a necessary and sufficient condition for a pair $(\alpha, \eta) \in \text{Aut}(H) \times \text{Aut}(K)$ to be inducible in a center extension.
	
	\begin{lemma}\label{First_lemma}
		Let $E(H,K)\equiv \xymatrix{e \ar[r] & H\ar[r]^{i} & G \ar[r]^{\beta} & K\ar[r] & e}$ be a center extension of multiplicative Lie algebras. Then the pair $(\alpha, \eta) \in Aut(H) \times Aut(K)$ is inducible if and only if there exist $\lambda: K \to H$ such that the following conditions hold
		\begin{enumerate}
			\item[$(1)$] $f(\eta(x),\eta(y))\alpha(f(x,y))^{-1} = \lambda(xy) \lambda(x)^{-1} \lambda(y)^{-1}$
			\item[$(2)$] $h(\eta(x), \eta(y)) \alpha(h(x,y))^{-1} = \lambda(x \star y ) \Gamma_{\eta(y)}(\lambda(x)) \Gamma_{\eta(x)}(\lambda(y)^{-1})$
			\item[$(3)$] $\Gamma_{\eta(x)}=\alpha \Gamma_x \alpha^{-1}$ or $ \alpha(t(x) \star h) = t(\eta(x)) \star \alpha(h)$ for all $x,y \in K$ and $h \in H$.
		\end{enumerate}
		\begin{proof}
			Let $\phi \in Aut_H(K)$ such that $\Pi(\phi)= (\alpha, \eta)$. Let $t:K \to G$ be a transversal. We know that any $g\in G$ can be written as $g=ht(x)$ for some $h\in H$ and $x\in K$. Consider $\beta(t(\eta(x))^{-1}\phi(t(x)))= \phi(x)^{-1} \phi(x)=e$. Since $Ker(\beta) = Img(i)=H$, this implies that there exist unique element $\lambda(x)$ (say) in $H$  such that $\lambda(x) = t(\eta(x))^{-1} \phi(t(x))$. Thus we get a map $\lambda: K \rightarrow H$. To prove $(1)$ and $(2)$, we use the fact that $\phi$ is a multiplicative Lie algebra homomorphism and $H$ is an abelian group with trivial multiplicative Lie product. 
			Let $g,g'\in G$. Then $g = ht(x)$ and $g^\prime=h^\prime t(y)$ for some $h,h'\in H$ and $x,y\in K$.\\
			(1). Since $\phi $ is a homomorphism, $\phi(gg')= \phi(g)\phi(g') $ this implies
			\begin{align*} 
				\phi(ht(x)h't(y)) &= \phi(ht(x))\phi(h't(y)) \\
				\phi(hh'f(x,y)t(xy)) & = \phi(ht(x))\phi(h't(y))\\
				\alpha(hh')\alpha(f(x,y)) \lambda(xy) t(\eta(xy)) & = \alpha(hh')\lambda(x)t(\eta(x))\lambda(y)t(\eta(y))\\
				\alpha(f(x,y))\lambda(xy) & = t(\eta(x)) t(\eta(y))  t(\eta(xy))^{-1} \lambda(x)\lambda(y)\\
				\alpha(f(x,y))\lambda(xy) & = t(\eta(x)) t(\eta(y))  t(\eta(xy))^{-1} \lambda(x)\lambda(y)\\
				\alpha(f(x,y))\lambda(xy) & = f(\eta(x),\eta(y)) \lambda(x)\lambda(y)\\
				f(\eta(x),\eta(y))\alpha(f(x,y))^{-1} & = \lambda(xy) \lambda(x)^{-1}\lambda(y)^{-1}.
			\end{align*}
			(2). Since $\phi(g\star g')=\phi(g)\star \phi(g')$, consider
			\begin{align*}
				\phi(ht(x) \star  kt(y)) & = \phi((h\star k)(t(x)\star k)(h\star t(y))(t(x)\star t(y)))\\
				& =\phi(t(x)\star k)\phi(h\star t(y))\phi(h(x,y)t(x\star y))\\
				& =(\phi(t(x)) \star \phi(k)) (\phi(h) \star \phi(t(y)))\alpha(h(x,y))\phi(t(x\star y))\\
				& =(\lambda(x) t(\eta(x))\star \alpha(k))(\lambda(y)t(\eta(y))\star \alpha(h))^{-1} \alpha(h(x,y))\lambda(x\star y)t(\eta(x\star y))\\
				& =(t(\eta(x))\star \alpha(k))(t(\eta(y))\star \alpha(h))^{-1}\alpha(h(x,y))\lambda(x\star y) t(\eta(x)\star \eta(y)).
			\end{align*}
			and
			\begin{align*}
				\phi(ht(x)) \star  \phi(kt(y)) & =\alpha(h)\phi(t(x))\star \alpha(k)\phi(t(y))\\
				& = \alpha(h)\lambda(x)t(\eta(x))\star \alpha(k)\lambda(y)t(\eta(y))\\
				& = (t(\eta(x))\star \alpha(k)\lambda(y)t(\eta(y)))(\alpha(h)\lambda(x)\star \alpha(k)\lambda(y) t(\eta(y)))\\
				& = (t(\eta(x))\star \alpha(k)\lambda(y))(t(\eta(x))\star t(\eta(y)))(\alpha(h) \lambda(x)\star t(\eta(y)))\\
				& =(t(\eta(x))\star \lambda(y))(t(\eta(x))\star \alpha(k))(t(\eta(x))\star t(\eta(y)))(\lambda(x)\star t(\eta(y)))\\
				& \phantom{=} ~~(\alpha(h)\star t(\eta(y))).		
			\end{align*}
			After comparing both
			\[h(\eta(x),\eta(y))\alpha(h(x,y))^{-1}=\lambda(x \star y)\Gamma_{\eta(y)}(\lambda(x)) \Gamma_{\eta(x)}(\lambda(y))^{-1}.\]
			(3). Let $x\in K$, $h \in H$. Then $\Gamma_{\eta(x)}(h)=t(\eta(x))\star h= \lambda(x) \phi(t(x))\star h= ~^{\lambda(x)}(\phi(t(x))\star h)(\lambda(x)\star h)=\phi(t(x)) \star \phi(\phi^{-1}(h))=\phi(t(x)\star \phi^{-1}(h))=\alpha(t(x)\star \alpha^{-1}(h))=\alpha(\Gamma_x(\alpha^{-1}(h)))=\alpha\Gamma_x\alpha^{-1}(h)$. Thus $\Gamma_\eta(x)=\alpha \Gamma_x \alpha^{-1}$. Or $\Gamma_{\eta(x)}\alpha=\alpha\Gamma_x$, this implies $ \alpha(t(x) \star h) = t(\eta(x)) \star \alpha(h)$. 
			%Consider	$ \phi(t(x) \star h) = \phi(t(x)) \star \phi(h) $, since $\alpha = \phi|_H$ implies $	\alpha(t(x) \star h)  = t (\eta(x)) \lambda(x) \star \alpha(h) = ~ ^{t(\eta(x))}(\lambda(x) \star \alpha(h))(t(\eta(x))\star \phi(h)) = (\lambda(x) \star \alpha(h)) (t(\eta(x)) \star \alpha(h)) ~ (\because H\subseteq Z(G)) $. Thus, we get $\alpha(t(x) \star h) = t(\eta(x)) \star \alpha(h)$.\\
			Conversely, since every $g\in G$ can be written as  $g=ht(x)$ for some $h\in H$ and $x\in K$. Then we can define a map $\phi:G \to G$ by $\phi(g) = \alpha(h)\lambda(x)t(\eta(x))$. Thus $\phi(h) = \alpha(h)$ for all $h \in H$ and $\eta(x)= \beta(t(\eta(x)))=\beta(\lambda(x)\phi(t(x)))=\beta(\phi(t(x)))$ for all $x\in K$. Now we prove that $\phi $ is an automorphism of G.  Let $g_1,g_2 \in G$. Then we can write $g_1=ht(x)$,  $g_2 = kt(y)$ for some $x,y\in K$ and $h,k\in H$. Consider
			\begin{align*}
				~~~~~~~~~\phi(ht(x)kt(y)) & = \phi(t(xy) f(x,y) hk)\\
				& = t(\eta(xy))\lambda (xy)\alpha(f(x,y)hk) \\
				& = t(\eta(x),\eta(y))\lambda(x)\lambda(y)\alpha(f(x,y))^{-1}f(\eta(x),\eta(y))\alpha(f(x,y)) \alpha(hk)\\
				& = t(\eta(x))t(\eta(y))f(\eta(x),\eta(y))^{-1}f(\eta(x),\eta(y))\lambda(x) \lambda(y)\alpha(h) \alpha(k)\\
				& = t(\eta(x))t(\eta(y))\lambda(x)\lambda(y)\alpha(h)\alpha(k)\\
				& = \phi(h t(x))\phi(kt(y)).
			\end{align*}
			For the operation $\star$ we have
			\begin{align*}
				\phi(t(x)h)\star \phi(t(y) k) & = t(\eta(x))\lambda(x)\alpha(h)\star t(\eta(y))\lambda(y)\alpha(k)\\
				& =~^{t(\eta(x))}(\lambda(x)\alpha(h)\star t(\eta(y))\lambda(y)\alpha(k))(t(\eta(x))\star t(\eta(y)) \lambda(y)\alpha(k))\\
				& = (\lambda(x)\alpha(h)\star t(\eta(y))\lambda(y)\alpha(k))(t(\eta(x))\star t(\eta(y))\lambda(y) \alpha(k))\\
				& =~^{\lambda(x)}(\alpha(h)\star t(\eta(y))\lambda(y)\alpha(k))(\lambda(x)\star t(\eta(y)) \lambda(y)\alpha(k))\\ & \phantom{=} ~~(t(\eta(x))\star t(\eta(y)))~^{t(\eta(y))}(t(\eta(x))\star \lambda(y)\alpha(k))\\
				& = (\alpha(h)\star t(\eta(y)))~^{t(\eta(y))}(\alpha(h)\star \lambda(y)\alpha(k))(\lambda(x)\star t(\eta(y)))\\
				& \phantom{=} ~~(t(\eta(x))\star t(\eta(y)))(t(\eta(x))\star\lambda(y)) ~ ^{\lambda(y)}(t(\eta(x))\star\alpha(k))\\
				& = (\alpha(h)\star t(\eta(y)))(\alpha(h)\star\lambda(y))(\alpha(h)\star \alpha(k))(\lambda(x)\star t(\eta(y)))\\
				& \phantom{=} ~~(\lambda(x)\star \lambda(y))(\lambda(x)\star  \alpha(k))(t(\eta(x))\star t(\eta(y)))(t(\eta(x))\star \lambda(y))\\
				& \phantom{=} ~~(t(\eta(x))\star \alpha(k))\\
				& = t(\eta(x \star y))h(\eta(x)\star \eta(y))(\lambda(x)\star t(\eta(y)))(t(\eta(x))\star \lambda(y))\\
				& \phantom{=} ~~(t(\eta(x))\star \alpha(k))(\alpha(h)\star t(\eta(y))).
			\end{align*}
			\begin{align*}
				\phi(t(x)h\star t(y)k) & = \phi((t(x)\star t(y))(t(x)\star k)(h\star t(y))(h\star k))\\
				& = \phi(t(x\star y)h(x,y))\alpha(t(x)\star k)\alpha(h\star t(y))\\
				& = t(\eta(x\star y))\lambda(x\star y)\alpha(h(x,y))(t(\eta(x))\star \alpha(k))(\alpha(h)\star t(\eta(y))) \\
				& = t(\eta(x\star y))h(\eta(x),\eta(y))(\lambda(x)\star t(\eta(y)))(t(\eta(x))\star \lambda(y))\\
				& \phantom{=}~ (t(\eta(x))\star \alpha(k))(\alpha(h)\star t(\eta(y))),
			\end{align*}
			Thus $\phi$ is a multiplicative Lie algebra homomorphism. Let $g \in ker(\phi)$, where $g=t(x)h$ for some $h\in H$ and $x \in K$. Then $\phi(g) = \phi(t(x)h) = t(\eta(x)) \lambda(x) \alpha(h)=e$, since $\lambda(x) \alpha(h) \in H$, this implies $t(\eta(x)) \in H$, i.e., $\beta(t(\eta(x))) =e$, $\eta(x)=e$, so $x=e$ as $\eta \in Aut(K)$, this implies $\alpha(h) =e$, i.e., $h=e$. Then  $g=e$ and $\phi $ is one-one.
			Now to  show $\phi $ is onto, let $g\in G$, $g = t(y) h'$, for some $y\in K$ and $h'\in H$. As $\eta \in Aut(K)$, there exists $y' \in K$ such that $\phi(y') = y$. Since $ \lambda(y')^{-1} h' \in H$, there exists $h \in H$ such that $\alpha (h) = \lambda(y')^{-1} h'$. After computing, we get $\phi(t(y')h) = t(\phi(y')) \lambda(y') \alpha(h) = t(y) \lambda(y') \lambda(y')^{-1} h'=t(y)h'=g$. Thus $\phi$ is onto. Therefore $\phi \in Aut_H(G)$.
	\end{proof}	\end{lemma}
	From the previous theorem, we can assert that every inducible pair $(\alpha, \eta) \in \text{Aut}(H) \times \text{Aut}(K)$ satisfies $\Gamma_{\eta(x)} = \alpha \Gamma_x \alpha^{-1}$, indicating that $\Gamma_{\eta(x)}$ is a homomorphism from $H$ to $H$. Consequently, we introduce the following definition.
	
	\begin{definition}\label{compatible_def}
		A pair $(\alpha, \eta) \in  Aut(H) \times Aut(K)$ is called compatible if $\alpha \Gamma_x \alpha^{-1} = \Gamma_{\eta(x)}$ for all $x \in K$. In other words, the following diagram commute
		\begin{center}
			\begin{tikzcd}
				K \arrow{r}{\eta} \arrow{d}{\Gamma} & K \arrow{d}{\Gamma}  \\
				End(H) \arrow{r}{\theta} & End(H)
			\end{tikzcd}
		\end{center}
		where $\theta(u) = \alpha u \alpha^{-1}~ \text{for all}~ u \in End(H)$. We denote the set of all compatible pair by $C^L$.
	\end{definition}
	\begin{remark}
		The set of all compatible pairs $C^{L} =\{(\alpha, \eta) \in Aut(H)\times Aut(K) : \alpha \Gamma_x \alpha^{-1} = \Gamma_{\eta(x)}~ \text{for all}~ x \in K\}$ is a subgroup of $Aut(H) \times Aut(K)$.
	\end{remark}
	\begin{remark}\label{remark}
		Since $\eta \in Aut(K)$, we can replace $x $ by $\eta^{-1}(x)$ in $(1)$, $(2)$ and get
		\begin{align*}
			f(x,y) \alpha(f(\eta^{-1}(x),\eta^{-1}(y)))^{-1} & = \lambda(\eta^{-1}(xy))\lambda(\eta^{-1}(x))^{-1}\lambda(\eta^{-1}(y))^{-1}\\
			f(x,y)(\alpha(f(\eta^{-1}(x),\eta^{-1}(y))))^{-1} & = \lambda'(xy)\lambda'(x)^{-1}\lambda'(y)^{-1}
		\end{align*}
		and
		\begin{align*}
			h(x,y)\alpha(h(\eta^{-1}(x),\eta^{-1}(y)))^{-1} & = \lambda(\eta^{-1}(x)\star \eta^{-1}(y)) \Gamma_y(\lambda(\eta^{-1}(x)))\Gamma_x(\lambda(\eta^{-1}(y)^{-1}))\\
			& = \lambda'(x\star y)\Gamma_y(\lambda'(x))\Gamma_x(\lambda'(y)^{-1}),
		\end{align*}
		where $\lambda'=\lambda\eta^{-1}$. It can be seen that $(\alpha(f(\eta^{-1},\eta^{-1})),\alpha(h(\eta^{-1},\eta^{-1})))\in Z_{ML(\Gamma)}^2(K,H)$. Thus multiplicative center $2$-cocycles $(f,h)$ and $(\alpha(f(\eta^{-1},\eta^{-1})),\alpha(h(\eta^{-1},\eta^{-1})))$ are cohomologous to each other. 
	\end{remark}
	\noindent From the Remark \ref{remark} and Lemma \ref{First_lemma} we have the following corollary.
	\begin{corollary}\label{first_corollary}
		Let $E(H,K)\equiv \xymatrix{e \ar[r]&H\ar[r]^{i}&G \ar[r]^{\beta} & K\ar[r]&e}$ be a center extension.  A compatible pair $(\alpha,\eta)$ is inducible if and only if   multiplicative center $2$-cocycles $(f,h)$ and $(\alpha\circ f \circ (\eta^{-1},\eta^{-1}),\alpha \circ h \circ(\eta^{-1},\eta^{-1}))$ are cohomologous.
	\end{corollary}
	An extension $E(H,K)\equiv \xymatrix{e \ar[r]&H\ar[r]^{i}&G \ar[r]^{\beta} & K\ar[r]&e}$ of multiplicative Lie algebra $H$ by $K$ is called a central extension if $H \subseteq Z(G) \cap LZ(G)$, where $LZ(G)$ is multiplicative center of $G$ and is defined as $LZ(G) = \{x \in G : x\star y =e ~ \text{for all} ~ y \in G\}$. Thus, if $E(H,K)$ is a central extension, then the map $\Gamma : K \to \text{End}(H)$ is $\Gamma_x(h) = t(x) \star h = e$ for all $x \in K$ and $h \in H$. For central extension, Lemma \ref{First_lemma} takes the following form.
	
	\begin{lemma}\label{third_lemma}
		Let $E(H,K)\equiv \xymatrix{e \ar[r] & H\ar[r]^{i} & G \ar[r]^{\beta} & K\ar[r] & e}$ be a central extension of multiplicative Lie algebras. Then the pair $(\alpha, \eta) \in Aut(H) \times Aut(K)$ is inducible if and only if there exist $\lambda: K \to H$ such that the following conditions hold
		\begin{enumerate}
			\item $f(\eta(x),\eta(y))\alpha(f(x,y))^{-1} = \lambda(xy)\lambda(x)^{-1}\lambda(y)^{-1}$
			\item $h(\eta(x), \eta(y)) \alpha(h(x,y))^{-1} = \lambda(x\star y)$ for all $x,y \in K$ and $h \in H$.
		\end{enumerate}
	\end{lemma}
	
	\section{Wells type exact sequence and  Inducible problem}
	In this section, we discuss the lifting of automorphisms in the context of center extensions of multiplicative Lie algebras. We establish the Wells type exact sequence and see that obstruction to inducibility of a pair of automorphisms in $Aut(H) \times Aut(K)$ lies in $H_{ML(\Gamma)}^2(K,H)$.
	
	Let $E(H,K)\equiv \xymatrix{e\ar[r] & H\ar[r]^{i} & G \ar[r]^{\beta} & K\ar[r] & e}$ be a center extension and $t: K \to G$ be a transversal of $E(H,K)$. Then $E(H,K)$ induced a multiplicative center $2$-cocycle $(f,h,\Gamma)\in Z^2_{ML(\Gamma)}(K,H)$. 
	Define $C_1^L=\{\alpha\in Aut(H):(\alpha,I)\in C^L\}$ and $C_2^L=\{\eta\in Aut(K):(I,\eta)\in C^L\}$. Thus we have
	\[\alpha \in C_1^L \iff \alpha \Gamma_x= \Gamma_x \alpha ~\text{or}~ \alpha(t(x) \star h) = t(x) \star \alpha(h),~ \text{for all}~ x \in K, ~h \in H\]
	and  \[\eta \in C_2^L \iff \Gamma_x=\Gamma_{\eta(x)}~ \text{or} ~t(x) \star h=t(\eta(x)) \star h,~\text{for all}~ x \in K,~h \in H.\]
	For any $\alpha \in C_1^L$ and $\eta \in C_2^L$. We can define maps $r_\alpha, r_\alpha',r_\eta, r_\eta': K \times K \to H$ by \[r_\alpha(x,y) = f(x,y) \alpha(f(x,y))^{-1}, ~r_\alpha'(x,y) = h(x,y) \alpha(h(x,y))^{-1}\] and
	\[r_\eta(x,y) = f(\eta(x),\eta(y))f(x,y)^{-1},~ r_\eta'(x,y) = h(\eta(x),\eta(y))h(x,y)^{-1}.\]
	
	\begin{lemma}\label{fifth_lemma}
		Let $(f,g,\Gamma)$ be a multiplicative center $2$-cocycle corresponding to the extension $E(H,K)$. Then $(r_\alpha,r_\alpha',\Gamma)$ and $(r_\eta,r_\eta',\Gamma)$ are multiplicative center $2$-cocycles, i.e. $(r_\alpha,r_\alpha',\Gamma)$ and $(r_\eta,r_\eta',\Gamma) \in Z_{ML(\Gamma)}^2(K,H)$.
		\begin{proof}
			To prove that  $(r_\alpha,r_\alpha',\Gamma) \in Z_{ML(\Gamma)}^2(K,H)$. We need to show that $(r_\alpha,r_\alpha',\Gamma)$ satisfies the equations from \ref{eq: 1} to \ref{eq: 7}. Clearly \ref{eq: 1} and \ref{eq: 2} hold. For \ref{eq: 3}, consider
			\begin{align*}
				r_\alpha(x,y) r_\alpha(xy,z) & = f(x,y)\alpha(f(x,y))^{-1}f(xy,z)\alpha(f(xy,z))^{-1}\\
				& = f(x,y)f(xy,z)\alpha(f(x,y)(xy,z))^{-1} \\
				& = f(y,z)f(x,yz)\alpha(f(y,z)f(x,yz))^{-1} \\
				& = f(y,z)\alpha(f(y,z))^{-1}f(x,yz)\alpha(f(x,yz))^{-1} \\
				& = r_\alpha(y,z)r_\alpha(x,yz).
			\end{align*}
			For \ref{eq: 4} consider \\
			$r_\alpha(y^{-1},y(x\star z))^{-1}r_\alpha(y(x\star z),y^{-1})r_\alpha(x\star y,y(x\star z)) r_\alpha'(x,y)r_\alpha'(x,z)$
			\begin{align*}
				& = f(y^{-1},y(x\star z))^{-1}\alpha(f(y^{-1},y(x\star z)))f(y(x\star z),y^{-1})\alpha(f(y(x\star z),y^{-1}))^{-1} \\
				& \phantom{=} ~~f(x\star y,~ ^y(x\star z))\alpha(f(x\star y ,~ ^y(x\star z)))^{-1}h(x,y) \alpha(h(x,y))^{-1} h(x,z)\alpha(h(x,z))^{-1}\\
				& = f(y^{-1},y(x\star z))^{-1}f(y(x\star z),y^{-1})f(x \star y,~ ^y(x\star z))h(x,y)h(x,z)  \alpha(f(y^{-1},y(x\star z)))\\
				& \phantom{=}~ \alpha(f(y(x\star z),y^{-1}))^{-1}\alpha(f(x\star y,~ ^y(x\star z)))^{-1}  \alpha(h(x,y))^{-1}\alpha(h(x,z))^{-1}\\
				& = \Gamma_x(f(y,z)) h(x,yz)\alpha(\Gamma_x(f(y,z))h(x,yz))^{-1}\\
				& = \Gamma_x(f(y,z)\alpha(f(y,z))^{-1})h(x,yz)\alpha(h(x,yz))^{-1}~(\text{For}~ \alpha \Gamma_x = \Gamma_x \alpha)\\
				& = \Gamma_x(r_\alpha(y,z))r_\alpha'(x,yz).
			\end{align*}
			In the same way, we can prove \ref{eq: 5}. Now for \ref{eq: 6} consider \\
			$\Gamma_x(r_\alpha(z^{-1},zy)^{-1}r_\alpha(zy,z^{-1})) \Gamma_y(r_\alpha(z^{-1},zx)r_\alpha(zx,z^{-1})^{-1})r_\alpha'(^zx,^zy)$
			\begin{align*}
				& = \Gamma_x(f(z^{-1},zy)^{-1}\alpha (f(z^-1,zy))f(zy,z^{-1})\alpha(f(zy,z^{-1}))^{-1}) \\
				& \phantom{=}~~ \Gamma_y(f(z^{-1},zx)\alpha(f(z^{-1},zx))^{-1} f(zx,z^{-1})^{-1} \alpha(f(zx,z^{-1})))h(^zx,^zy)\alpha(h(^zx,^zy))^{-1} \\
				& = \Gamma_x(f(z^{-1},zy)^{-1}f(zy,z^{-1}))\Gamma_y(f(z^{-1},zx)f(zx,z^{-1})^{-1})h(^zx,^zy)\\
				& \phantom{=}~~ \alpha(\Gamma_x(f(z^{-1},zy)f(zy,z^{-1})^{-1}) \Gamma_y(f(z^{-1},zx)^{-1})f(zx,z^{-1})h(^zx,^zy)^{-1})\\
				& = f(z^{-1},z(x\star y))^{-1} f(z(x\star y),z^{-1})h(x,y)\alpha(f(z^{-1},z(x\star y))f(z(x\star y),z^{-1})^{-1} h(x,y)^{-1})\\
				& = r_\alpha(z^{-1},z(x\star y))^{-1}r_\alpha(z(x\star y),z^{-1})r_\alpha'(x,y).
			\end{align*}
			By the similar computation relation \ref{eq: 7} holds. Thus we get $(r_\alpha, r_\alpha', \Gamma) \in Z_{ML(\Gamma)}^2(K,H)$. Similarly, we can prove that $(r_\eta,r_\eta',\Gamma) \in Z_{ML(\Gamma)}^2(K,H)$.
		\end{proof}
	\end{lemma}
	By Lemma \ref{fifth_lemma}, we can define the following maps $\chi_1 : C_1^L \rightarrow H_{ML(\Gamma)}^2(K,H)$ and $\chi_2 : C_2^L \to H_{ML(\Gamma)}^2(K,H)$ by 
		\begin{align*}
		\chi_1(\alpha)= & (r_\alpha,r_\alpha',\Gamma)B_{ML(\Gamma)}^2(K,H) \\
	\chi_2(\eta)= & (r_\eta,r_\eta',\Gamma)B_{ML(\Gamma)}^2(K,H).
	\end{align*}

	\begin{lemma}\label{sixth_lemma}
		If $(r_\alpha,r_\alpha',\Gamma)$ and $(r_\eta,r_\eta',\Gamma)$ are multiplicative center $2$-cocycle. Then the maps $\chi_1$ and $\chi_2$ are well defined.
		\begin{proof}
			To prove that $\chi_1$ and $\chi_2$ are well defined, we need to show that $\chi_1$ and $\chi_2$ are independent of the transversal $t$. Let $t$ and $t'$ be two transversals of the extension $E(H,K)\equiv \xymatrix{e \ar[r]&H\ar[r]^{i}&G \ar[r]^{\beta}&K\ar[r] & e}$ with $t(e)=e$ and $t'(e)=e$. Then there exists maps $f,f',h~ \text{and}~h': K \times K \to H$ defined by $f(x,y) = t(x)t(y)t(xy)^{-1}, f'(x,y)=t'(x)t'(y)t'(xy)^{-1},h(x,y) = (t(x)\star t(y)) t(x \star y)^{-1}$ and $h'(x,y) = (t'(x) \star t'(y)) t'(x \star y)^{-1}$. Since $\beta(t(x)) = \beta(t'(x))=x$, this implies $t(x)(t'(x))^{-1} \in Ker(\beta) = Img(i)= H$. Thus there exist a unique $\lambda(x) \in H$, i.e. we get $\lambda: K \to H$ such that $t(x) = \lambda(x) t'(x)$ for all $x \in K$. Using $\lambda$ we get,
			\begin{align*}
					f(x,y)& = \lambda(x) \lambda(y) \lambda(xy)^{-1} f'(x,y) ~ \text{and}\\
				h(x,y)& = \Gamma_x(\lambda(y)) \Gamma_y(\lambda(x))^{-1} \lambda(x \star y)^{-1} h'(x,y).
			\end{align*}
			Let $\alpha \in C_1^L$ and $\eta \in C_2^L$. Suppose maps $r_\alpha,r_\alpha',r_\eta,r_\eta' : K \times K \to H$ associated with $t$ are same as defined before. For transversal $t'$ we define $s_\alpha,s_\alpha',s_\eta,s_\eta' : K \times K \to H$ as $s_\alpha(x,y) = f'(x,y) \alpha(f'(x,y))^{-1}, ~s_\alpha'(x,y) = h'(x,y) \alpha(h'(x,y))^{-1}, ~s_\eta(x,y) = f'(\eta(x),\eta(y))f'(x,y)^{-1}$ and $ s_\eta'(x,y) = h'(\eta(x),\eta(y))h'(x,y)^{-1}.$ Consider
			\begin{align*}
				r_\alpha(x,y) s_\alpha(x,y)^{-1} & = f(x,y) \alpha(f(x,y))^{-1}(f'(x,y) \alpha(f'(x,y))^{-1})^{-1} \\
				& = 	h(x,y) h'(x,y)^{-1} \alpha(h(x,y) h'(x,y)^{-1} )^{-1} \\
				& = (\lambda(x) \lambda(y) \lambda(xy)^{-1}) \alpha(\lambda(x) \lambda(y) \lambda(xy)^{-1})^{-1}
			\end{align*}
			On the other hand
			\begin{align*}
				r'_\alpha(x,y) s'_\alpha(x,y)^{-1} & = h(x,y) \alpha(h(x,y))^{-1}(h'(x,y) \alpha(h'(x,y))^{-1})^{-1} \\
				& = h(x,y) h'(x,y)^{-1} \alpha(h(x,y) h'(x,y)^{-1} )^{-1}\\
				& = \Gamma_x(\lambda(y)) \Gamma_y(\lambda(x))^{-1} \lambda(x \star y)^{-1} \alpha(\Gamma_x(\lambda(y)) \Gamma_y(\lambda(x))^{-1} \lambda(x \star y)^{-1})^{-1}\\
				& = \Gamma_x(\lambda(y) \alpha(\lambda(y))^{-1}) \Gamma_y(\lambda(x) \alpha(\lambda(x))^{-1})^{-1} \lambda(x \star y)^{-1} (\alpha(\lambda(x \star y))^{-1})^{-1}.
			\end{align*}
			Define a map $\lambda': K \to H$ by $\lambda'(x) = \lambda(x) \alpha (\lambda(x))^{-1}$, clearly $\lambda'(e)=e $. Thus we get $r_\alpha(x,y) s_\alpha(x,y)^{-1} = \lambda'(x) \lambda'(y) \lambda'(xy)^{-1}$ and $r'_\alpha(x,y) s'_\alpha(x,y)^{-1}= \Gamma_x(\lambda'(y)) \Gamma_y(\lambda'(x))^{-1}\\ \lambda'(x \star y)^{-1}$. Thus $(r_\alpha s_\alpha^{-1}, r_\alpha's_\alpha'^{-1}, \Gamma) \in B_{ML(\Gamma)}^2(K,H)$. Hence the map $\chi_1$ is independent of transversal. Similarly we can show that $\chi_2$ is independent of transversal. Therefore both $\chi_1$ and $\chi_2$ are well defined.
		\end{proof}
	\end{lemma}
	\noindent \textbf{Note:} The maps $\chi_1$ and $\chi_2$ need not be homomorphisms as group operation for $C_1^L$ and $C_2^L$ is composition not pointwise addition. Whereas $Ker(\chi_1)$ and $Ker(\chi_2)$ have usual meaning.\\
	\textbf{Notations:} Here we define some notations which we use in the entire paper.
	\begin{align*}
		Aut_H(G) & = \{ \phi \in Aut(G) : \phi(H) = H\}\\
		Aut^H(G) & = \{ \phi \in Aut(G): \phi(h)= h ~ \text{for all}~ h \in H\}\\
		Aut^{H,K}(G) & = \{ \phi \in Aut(G) : \phi(h) =h ~ \text{and} ~ \phi(x)x^{-1} \in H ~ \text{for all} ~ x \in G, ~ h \in H\} \\
		Aut^H_K(G) & = \{ \phi \in Aut(G) : \phi(H) =H ~ \text{and} ~ \phi(x)x^{-1} \in H ~ \text{for all} ~ x \in G\}.
	\end{align*}

	Now with the help of above discussion. We are set to address the following problem in the next theorem for multiplicative Lie algebra.\\
	\textbf{Problem:} Let $E(H,K)\equiv \xymatrix{e \ar[r]&H\ar[r]^{i}&G \ar[r]^{\beta}&K\ar[r] & e}$ be a center extension. Under what conditions \begin{enumerate}
		\item[(1)] a pair $(\alpha,I) \in C^L$ is inducible. Or, $\alpha \in C_1^L$ can be lifted to an automorphism of $G$, which fixes $H$.
		\item[(2)] a pair $(\eta,I) \in C^L$ is inducible. Or, $\eta \in C_2^L$ can be lifted to an automorphism $\phi$ of $G$ such that $\beta(\phi(t(x)))=x$ for all $x\in K$.
	\end{enumerate}
	\begin{theorem}\label{first_thm}
		Let $E(H,K)\equiv \xymatrix{e \ar[r]&H\ar[r]^{i}&G \ar[r]^{\beta}&K\ar[r] & e}$ be a center extension. Then there exist the following two exact sequences
		\begin{eqnarray*}
			e \parbox{0.8cm}{\rightarrowfill} Aut^{H,K}(G)  \xrightarrow{\makebox[0.8cm]{$i$}} Aut_H^K(G) \xrightarrow{\makebox[0.8cm]{$\Pi_1$}}   C_1^L \xrightarrow{\makebox[0.8cm]{$\chi_1$}} H^2_{ML(\Gamma)}(K,H)\\
			e \parbox{0.8cm}{\rightarrowfill} Aut^{H,K}(G)  \xrightarrow{\makebox[0.8cm]{$i$}} Aut^H(G) \xrightarrow{\makebox[0.8cm]{$\Pi_2$}}  C_2^L \xrightarrow{\makebox[0.8cm]{$\chi_2$}} H^2_{ML(\Gamma)}(K,H),
		\end{eqnarray*}
		where the maps $\Pi_1(\phi)= \alpha$ and $\Pi_2(\phi) = \eta$.
		\begin{proof}
			Clearly, both the sequences are exact at the first two terms. Now we prove the exactness of the first sequence at the third term. Let $\phi \in Aut_H^K(G)$. Then $\Pi_1(\phi) \in C_1^L$ and $(\chi_1\Pi_1)(\phi) = \chi_1(\Pi_1(\phi)) = \chi_1(\alpha) = (r_\alpha, r_\alpha', \Gamma) B_{ML(\Gamma)}^2(K,H)$, where $\alpha = \phi|_H$. Since $\phi \in Aut_H^K(G)$, we have $\eta(x) = (\beta \phi t)(x) = x $. Now we have
			\begin{align*}
				r_\alpha(x,y ) & = f(x,y) \alpha(f(x,y))^{-1}\\
				& = f(\eta(x),\eta(y)) \alpha (f(x,y))^{-1}~ (\because~ \eta(x) =x)\\
				& = \lambda(xy) \lambda(x)^{-1} \lambda(y)^{-1} ~ (\text{by Lemma}~ \ref{First_lemma})
			\end{align*}
			\begin{align*}
				r_\alpha'(x,y ) & = h(x,y) \alpha(h(x,y))^{-1}\\
				& = h(\eta(x),\eta(y)) \alpha (h(x,y))^{-1}~ (\because~ \eta(x) =x)\\
				& = \lambda(x \star y) \Gamma_y(\lambda(x)) \Gamma_x(\lambda(y))^{-1} ~ (\text{by Lemma}~ \ref{First_lemma}).
			\end{align*}
			Thus $(r_\alpha,r_\alpha',\Gamma) \in B_{ML(\Gamma)}^2(K,H,\Gamma)$, this implies $Img(\Pi_1) \subseteq ker(\chi_1)$. Let $\alpha \in C_1^L$. Then $(\alpha,I)\in C^L$ such that $\chi_1(\alpha) \in B_{ML(\Gamma)}^2(K,H)$. Remark \ref{remark} implies that cocycles $(f,h)$ and $(\alpha \circ f, \alpha \circ h)$ are cohomologous. Thus Corollary \ref{first_corollary} implies that $(\alpha,I)$ is inducible. Then there exists $\phi\in Aut_H(G)$ such that $\Pi(\phi)=(\alpha,I)$. Clearly $\phi\in Aut_H^K(G)$ with $\Pi_1(\phi)=\alpha$. Thus $ker(\chi_1) \subseteq Img(\phi_1)$.  
			%there exist a map $\lambda : K \to H$, $\lambda(e) =e $ such that $r_\alpha(x,y) = \lambda(xy)^{-1} \lambda(x) \lambda(y)$ and $r_\alpha'(x,y) = \Gamma_x(\lambda(y)) \Gamma_y (\lambda(x))^{-1} \lambda(x \star y)$, Set $\eta = id|_K$, then by Lemma \ref{First_lemma} there exist $\phi: G \to G$ such that $\phi(ht(x)) = t(\eta(x)) \lambda(x) \alpha(h) = t(x) \lambda(x) \alpha(h)$ for all $x \in K$ and $h \in H$. Compute $\phi(t(x) h) t(x)^{-1} h^{-1} =t(x) \lambda(x) \alpha(h) t(x)^{-1} h^{-1} = \lambda(x) \alpha(h) h^{-1} \in H $, this implies $\phi(g) g^{-1} \in H$, i.e., $\phi \in Aut_H^K(G)$. Also $\phi(h) = \alpha(h)$ for all $h \in H$, this implies $\Pi_1(\phi) = \alpha $ and $\alpha\in Img(\Pi_1)$. Thus $ker(\chi_1) \subseteq Img(\Pi_1)$.   \\
			Now we prove the exactness for the second sequence in the third term. Let $\phi \in Aut^H(G)$, $\Pi_2(\phi) = \eta \in C_2^L$  and $(\chi_2 \Pi_2)(\phi) = \chi_2(\Pi_2(\phi)) = \chi_2(\eta)= (r_\eta,r_\eta',\Gamma)B_{(ML(\Gamma))}^2(K,H)$, where $\eta(x) = (\beta \phi t)(x)$ for all $x \in K$. Since $ \eta = \Pi_2( \phi)\in C_2^L$, $\Gamma_x = \Gamma_{\eta(x)}$ for all $x \in K$, $\alpha= \phi|_H = id_H$ for $\phi \in Aut^H(G)$. 
			\begin{align*}
				r_\eta(x,y) & = f(\eta(x),\eta(y)) f(x,y)^{-1}\\
				& = f(\eta(x),\eta(y)) \alpha(f(x,y))^{-1} \\
				& = \lambda(xy) \lambda(x)^{-1} \lambda(y)^{-1} ~ (\text{by Lemma}~ \ref{First_lemma})
			\end{align*}
			\begin{align*}
				r_\eta'(x,y) & = h(\eta(x),\eta(y)) h(x,y)^{-1}\\
				& = h(\eta(x) , \eta(y)) \alpha(h(x,y))^{-1} \\
				& = \lambda(x \star y ) (\lambda \star t(\eta(y)))^{-1} (t(\eta (x)) \star \lambda(y))^{-1} ~ (\text{by Lemma}~ \ref{First_lemma})\\
				& = \lambda(x \star y) \Gamma_y (\lambda(x)) \Gamma_x(\lambda(y))^{-1}
			\end{align*}
			This implies $(r_\eta, r_\eta', \Gamma)\in B_{ML(\Gamma)}^2(K,H)$. Thus $Img(\Pi_2) \subseteq ker(\chi_2)$. Let $\eta \in ker(\chi_2)$. Then $(I,\eta)\in C^L$ such that $\chi_2(\eta) \in B_{ML(\Gamma)}^2(K,H)$. Remark \ref{remark}, implies that cocycles $(f,h)$ and $(f\circ(\eta^{-1},\eta^{-1}),h\circ(\eta^{-1},\eta^{-1}))$ are cohomologus. Thus Corollary \ref{first_corollary} implies that $(I,\eta)$ is inducible. Then there exist $\phi\in Aut_H(G)$ such that $\Pi(\phi)=(I,\eta)$. Clearly $\eta \in Aut_H^K(G)$ with $\Pi_2(\phi)=\eta$. Thus $Ker(\chi_2)\subseteq Img(\phi_2)$.
			% there exist $\lambda : K \to H$, $\lambda(e) =e$ such that $r_\eta(x,y) = \lambda(x) \lambda(y) \lambda(xy)^{-1}$ and $r_\eta'(x,y) = \Gamma_x(\lambda(y)) \Gamma_y(\lambda(x))^{-1} \lambda(x \star y)^{-1}$. Let $\alpha = id|_H$, then by Lemma \ref{First_lemma} there exist an automorphism $\phi : K \to H$ such that $\phi(t(x) h) = t (\eta(x)) \lambda(x) \alpha(h)$ for all $ x\in k$ and $ h \in H$. so $\phi ( h) = \phi(t(e) h) = \alpha(h) =h $, implies $\phi \in Aut^H(G)$ and $\Pi_2(\phi ) = \eta $. Thus $ker(\chi) \subseteq Img(\Pi_2)$.
		\end{proof}
	\end{theorem}
	\begin{remark}
		From the last Theorem \ref{first_thm} we conclude that $\alpha \in C_1^L$ and $\eta \in C_2^L$ can be lifted to an automorphism in $Aut_H(G)$, respectively if and only if $\alpha \in Ker(\chi_1)$ and $\eta \in Ker(\chi_2)$. 
	\end{remark}
	
	Let $E(H,K)\equiv \xymatrix{e \ar[r]&H\ar[r]^{i}&G \ar[r]^{\beta}&K\ar[r] & e}$ be a center extension and $(\alpha, \eta) \in C^L $. Define maps $s_{\alpha, \eta} , s_{\alpha, \eta}' : K \times K \to H$  by 
	\[s_{\alpha,\eta }(x,y) = r_\alpha(x,y) r_\eta(x,y) = f(\eta(x),\eta(y)) \alpha(f(x,y))^{-1}\] and \[s_{\alpha, \eta}'(x,y) = r_\alpha'(x,y) r_\eta'(x,y) = h(\eta(x),\eta(y)) \alpha(h(x,y))^{-1}.\]
	% Since $Z_{ML(\Gamma)}^2(K,H)$ is an abelian group. Thus by 
	Using Lemma \ref{fifth_lemma}, $(s_{\alpha,\eta },s_{\alpha, \eta}', \Gamma) \in Z_{ML(\Gamma)}^2(K,H)$. We know that different $2$-cocycles associated with the same extension differ by an $2$-coboundry. Thus we can define a map $\chi: C^L\to H^2_{ML(\Gamma)}(K,H)$ by \[\chi(\alpha,\eta) = (s_{\alpha,\eta },s_{\alpha, \eta}',\Gamma) B_{ML(\Gamma)}^2(K,H). \]
	From Lemma \ref{sixth_lemma} we observe that $\chi$ is independent of the choice of transversal $t$. Thus $\chi$ is well defined. 
	\begin{remark}
		If $\chi(\alpha,\eta)$ is trivial for the pair $(\alpha,\eta)\in C^L$, where $\eta \in \text{Aut}(K)$, we can think $x$ as $\eta^{-1}(x)$. From Corollary \ref{first_corollary}, $(f,h)$ and $(\alpha\circ f \circ (\eta^{-1},\eta^{-1}),\alpha \circ h \circ(\eta^{-1},\eta^{-1}))$ are cohomologous cocycles. Thus, the compatible pair $(\alpha,\eta)\in C^L$ is inducible if and only if $\chi(\alpha,\eta)$ is trivial. Therefore, we can assert that $\chi(\alpha,\eta)$ is an obstruction for the inducibility of the compatible pair $(\alpha,\eta)$.
	\end{remark}
	
	\begin{lemma}
		Let $E(H,K)\equiv \xymatrix{e \ar[r]&H\ar[r]^{i}&G \ar[r]^{\beta}&K\ar[r] & e}$ be a center extension. Then $Z_{ML(\Gamma)}^1(K,H) \cong Aut^{H,K}(G)$ as groups.
		\begin{proof}
			Define $\delta: Z_{ML(\Gamma)}^1(K,H) \to Aut(G)$ by $\delta(\nu) = \phi_\nu$, where $\phi_\nu: G \to G$ is defined by $\phi_\nu(ht(x))=h\nu(x) t(x)$, $h\in H$ and $x\in K$. To show $\phi_\nu \in Aut(G)$. Let $g=ht(x)$ and $g'=h't(y)\in G$. First consider
			\begin{align*}
				\phi_\nu(ht(x) h't(y)) & = \phi_\nu(hh't(x)t(y)) \\
				& = \phi_\nu(hh'f(x,y)t(xy))\\
				& = hh'f(x,y) \nu(xy)t(xy) \\
				& = hh'f(x,y) \nu(x)\nu(y) f(x,y)^{-1}t(x)t(y)\\
				& = h \nu(x) t(x) h'\nu(y) t(y)\\
				& = \phi_\nu(ht(x)) \phi_\nu(h't(y)).
			\end{align*}
			Now for the operation $\star$ we have
			\begin{align*}
				\phi_\nu(ht(x) \star h't(y)) & = \phi_\nu((h\star h')(h\star t(y))(t(x)\star h')(t(x)\star t(y))) \\
				& = \phi_\nu((h\star t(y))(t(x)\star h')h(x,y)(t(x\star y)))\\
				& = (t(y)\star h)^{-1} (t(x)\star h')h(x,y) \nu(x\star y) t(x \star y) \\
				& = (t(y)\star h)^{-1} (t(x)\star h')h(x,y) (t(x)\star \nu(y)) (t(y) \star \nu(x))^{-1}h(x,y)^{-1}\\
				& \phantom{=} ~~(t(x) \star t(y))\\
				& = (t(y)\star h)^{-1} (t(x)\star h') (t(x)\star \nu(y)) (t(y) \star \nu(x))^{-1} (t(x) \star t(y))\\
				& = (h\star h') (h \star t(y)) (t(x)\star h') (t(x)\star \nu(y)) (\nu(x) \star t(y)) (t(x) \star t(y))\\
				& = (h\nu(x) t(x) \star h'\nu(y) t(y)) \\
				& = 	\phi_\nu(ht(x)) \star \phi_\nu(h't(y)).
			\end{align*}
			Thus $\phi_\nu$ is a multiplicative Lie algebra homomorphism. Let $g=ht(x) \in G$ and $\phi_\nu(ht(x))=h\lambda(x)t(x)=e$. Then $t(x)\in H$, $x=\beta(t(x))=e$ and $h=e$, $g=e$. Thus $\phi_\nu$ is injective. If $ht(x) \in G$, then $\phi_\nu(h\lambda(x)^{-1}t(x))=ht(x)$. Thus $\phi_\nu$ is surjective. Hence $\phi_\nu \in Aut(G)$. Clearly, $\phi_\nu(h) =h$ for all $h\in H$ and $\phi_\nu(ht(x))(ht(x))^{-1}=h\nu(x)t(x)t(x)^{-1}h^{-1}=\nu(x) \in H$. Therefore $\phi_\nu \in Aut^{H,K}(G)$.\\
			Now we prove that $\delta$ is a group isomorphism. Let $\nu_1$ and $\nu_2\in Z_{ML(\Gamma)}^1(K,H)$. Then
			\begin{align*}
				\delta(\nu_1+\nu_2)(g) & = \phi_{\nu_1+\nu_2}(ht(x)) \\
				& =h(\nu_1+\nu_2)(x)t(x) \\
				& = h\nu_1(x) t(x) + h\nu_2(x)t(x) \\
				& = (\delta(\nu_1)+\delta(\nu_2))(g).
			\end{align*}
			Thus $\delta$ is a group homomorphism. Let $g\in G$ such that $\delta(\nu)(g)=I_\nu(g)=g$, this implies $h\lambda(x) t(x)=ht(x)$, $\lambda(x)=e$ for all $x\in K$. Thus $\delta$ is injective. Let $\phi \in Aut^{H,K}(G)$. Then $\beta(\phi(t(x))t(x)^{-1})=e$, this implies $\phi(t(x))t(x)^{-1}\in Ker(\beta)=H$, there exists unique element say $\nu(x) \in H$ such that $\phi(t(x))=\nu(x)t(x)$. Thus we have a map $\nu:K\to H$ such that $\phi(t(x))=\nu(x)t(x)$. To show that $\nu \in Z_{ML(\Gamma)}^1(K,H)$, we use the fact that $\phi$ is a multiplicative Lie algebra homomorphism. Consider
			\begin{align*}
				\phi(t(x)t(y))& =\phi(t(x))\phi(t(y))\\
				\phi(f(x,y)t(xy))	& =\nu(x) t(x) \nu(y)t(y) \\
				f(x,y) \nu(xy) t(xy)	& = \nu(x) \nu(y) f(x,y) t(xy) \\
				\nu(xy) 	& = \nu(x) \nu(y).
			\end{align*}
			Thus $\nu$ is a group homomorphism. Now consider the map $\phi$ for the operation $\star$
			\begin{align*}
				\phi(t(x)\star t(y))& =\phi(t(x))\star \phi(t(y))\\
				\phi(h(x,y) t(x\star y))	& =\nu(x) t(x)\star \nu(y)t(y) \\
				h(x,y) \nu(x\star y) t(x\star y)	& = (\nu(x)\star \nu(y)) (\nu(x)\star t(y))(t(x)\star \nu(y)) (t(x)\star t(y)) \\
				\nu(x\star y) 	& =(t(x)\star \nu(y)) (t(y)\star \nu(x))^{-1} \\
				\nu(x\star y) 	& =\Gamma_x(\nu(y)) \Gamma_y(\nu(x)^{-1}).
			\end{align*}
			Also $\phi(ht(x)) = \phi(h) \phi(t(x))= h \nu(x) t(x)$, $h\in H$ and $x\in K$. Thus $\nu\in Z_{ML(\Gamma)}^1(K,H)$. Hence $\delta$ is surjective for $Aut^{H,K}(G)$. Therefore $Z_{ML(\Gamma)}^1(K,H) \cong Aut^{H,K}(G)$.
		\end{proof}
	\end{lemma}
	
	\begin{theorem}\label{second_thm}
		Let $E(H,K)\equiv \xymatrix{e \ar[r]&H\ar[r]^{i}&G \ar[r]^{\beta}&K\ar[r] & e}$ be a center extension. Then we have the following exact sequence 
		\[ e \parbox{0.8cm}{\rightarrowfill} Z_{ML(\Gamma)}^1(K,H)   \xrightarrow{\makebox[0.8cm]{$i$}} Aut_H(G) \xrightarrow{\makebox[0.8cm]{$\Pi$}}   C^L \xrightarrow{\makebox[0.8cm]{$\chi$}} H^2_{ML(\Gamma)}(K,H),~~~~~~~~~~(A)\] where $\chi(\alpha, \eta) = (s_{\alpha,\eta },s_{\alpha, \eta}',\Gamma) B_{ML(\Gamma)}^2(K,H)$ and $ \Pi(\phi) =(\alpha, \eta)$.
		\begin{proof} Clearly, the sequence is exact at the first two terms. For third term, let $\phi \in Aut_H(G)$ such that $\Pi(\phi) =(\alpha, \eta)$. Then by Lemma \ref{First_lemma}, $f(\eta(x),\eta(y))\alpha(f(x,y))^{-1} = \lambda(xy) \lambda(x)^{-1} \lambda(y)^{-1}$ and $h(\eta(x),\eta(y)) \alpha(h(x,y))^{-1} = \lambda(x \star y ) \Gamma_{\eta(y)}(\lambda(x)) \Gamma_{\eta(x)}(\lambda(y)^{-1})$, where $\lambda:K \to H$, this implies $s_{\alpha,\eta }(x,y) = \lambda(xy) \lambda(x)^{-1} \lambda(y)^{-1}$ and $s_{\alpha, \eta}'(x,y) = \lambda(x \star y ) \Gamma_{\eta(y)}(\lambda(x)) \Gamma_{\eta(x)}(\lambda(y)^{-1})$. Thus $(s_{\alpha,\eta},s_{\alpha,\eta}', \Gamma) \in B_{ML(\Gamma)}^2(K,H)$. Hence $(\alpha,\eta) \in ker(\chi)$.
			Let $(\alpha,\eta) \in ker(\chi)$ such that $\chi(\alpha,\eta) = (s_{\alpha,\eta },s_{\alpha, \eta}',\Gamma) B_{ML(\Gamma)}^2(K,H) =  B_{ML(\Gamma)}^2(K,H)$, this implies $(s_{\alpha,\eta},s_{\alpha,\eta}',\Gamma) \in B_{ML(\Gamma)}^2(K,H)$, then there exist $\lambda: K\to H$ with $\lambda(e)=e$ such that $s_{\alpha,\eta}(x,y) = \lambda(xy)^{-1} \lambda(x) \lambda(y)$ and $s_{\alpha,\eta}'(x,y) = \lambda(x \star y ) \Gamma_{\eta(y)}(\lambda(x)) \Gamma_{\eta(x)}(\lambda(y)^{-1})$. Then by Lemma \ref{First_lemma}, the pair $(\alpha,\eta)$ is inducible. Thus there exist $\phi \in Aut_H(G)$ such that $\Pi(\phi) = (\alpha, \eta)$. Hence $(\alpha ,\eta ) \in Img(\Pi)$.
		\end{proof}
	\end{theorem}
	
	\begin{corollary}
		Let $E(H,K)\equiv \xymatrix{e \ar[r]&H\ar[r]^{i}&G \ar[r]^{\beta}&K\ar[r] & e}$ be a split center extension. Then every compatible pair $(\alpha,\eta) \in Aut(H) \times Aut(K)$ is inducible.
		\begin{proof}
			Since $E(H,K)$ is a split extension. This implies that maps $f: K \times K \to H$ and $g: K \times K \to H$ both are trivial, then the map $\chi $ is also trivial. From the exactness of the sequence $(A)$, we get $\Pi$ is surjective.
		\end{proof}
	\end{corollary}
	
	\begin{theorem}\label{second_second_thm}
		Let $E(H,K)\equiv \xymatrix{e \ar[r]&H\ar[r]^{i}&G \ar[r]^{\beta}&K\ar[r] & e}$ be a central extension. Then we have the following exact sequence 
		\[ e \parbox{0.8cm}{\rightarrowfill} Z_{ML(\Gamma)}^1(K,H)  \xrightarrow{\makebox[0.8cm]{$i$}} Aut_H(G) \xrightarrow{\makebox[0.8cm]{$\Pi$}}   Aut(H) \times Aut(K) \xrightarrow{\makebox[0.8cm]{$\chi$}} H^2_{ML(\Gamma)}(K,H),\] where $\chi(\alpha, \eta) = (s_{\alpha,\eta },s_{\alpha, \eta}',\Gamma) B_{ML(\Gamma)}^2(K,H)$ and $ \Pi(\phi) =(\alpha, \eta)$.
		\begin{proof}
			Suppose $E(H,K)\equiv \xymatrix{e \ar[r]&H\ar[r]^{i}&G \ar[r]^{\beta}&K\ar[r] & e}$ is a central extension, i.e. $H \subseteq Z(G) \cap LZ(G)$. Thus we have $\Gamma_x(h) = t(x) \star h = e$, this implies that if $(\alpha,\eta) \in Aut(H)\times Aut(K)$, we have $\alpha \Gamma_x \alpha^{-1} = \Gamma_{\eta(x)}$, i.e. $(\alpha,\eta) \in C^L$. Hence for a central extension $Aut(H)\times Aut(K) = C^L$. Similar to the Theorem \ref{second_thm}, the exactness of the sequence follows using Lemma \ref{third_lemma}.
		\end{proof}
	\end{theorem}
	
	\begin{corollary}
		Let $E(H,K)\equiv \xymatrix{e \ar[r]&H\ar[r]^{i}&G \ar[r]^{\beta}&K\ar[r] & e}$ be a split central extension. Then every pair $(\alpha,\eta) \in Aut(H) \times Aut(K)$ is inducible.
	\end{corollary}
	\begin{definition} \cite{RLS} Let 
		$E(R,K) \equiv \xymatrix{e \ar[r] & R\ar[r]^{i} & F \ar[r]^{\nu} & K\ar[r] & e}$
		be a free presentation of a multiplicative Lie algebra $K$. Then Schur Multiplier of the multiplicative Lie algebra $K$ is defined by $\frac{([F,F](F \star F))\cap R}{[R,F](R \star F)}$ and denoted by $\tilde{M}(K)$.
	\end{definition}
	Recall from [\cite{RLS}, Corollary 6.3], If $K$ is a finite multiplicative Lie algebra, then $\tilde{M}(K) \cong H_{ML}^2(K, \mathbb{C}^\star)$. Also, a multiplicative Lie algebra $K$ is said to be perfect if $K=[K,K](K\star K)$. Now we mentioned two Propositions from \cite{RLS} as we will use them in our next corollary.
	\begin{proposition}\label{P_1} \cite{RLS} Every perfect multiplicative Lie algebra $K$ admits a unique universal central extension. More precisely, If free presentation of $K$ is given by the extension $E(R,K)$ then universal central extension by $K$ is given by
		\[ 	e \parbox{0.8cm}{\rightarrowfill} \frac{([F,F](F \star F))\cap R}{[R,F](R \star F)}  \xrightarrow{\makebox[0.8cm]{$\bar{i}$}} \frac{([F,F](F \star F))}{[R,F](R \star F)} \xrightarrow{\makebox[0.8cm]{$\bar{\nu}$}}   K \longrightarrow  e. \]
	\end{proposition}
	\begin{proposition}\label{P_2} \cite{RLS} A central extension $E(H,K) \equiv \xymatrix{e \ar[r] & H\ar[r]^{i} & G \ar[r]^{\beta} & K\ar[r] & e}$ is a universal central extension by $K$ if and only if $G$ is perfect and every central extension by $G$ splits.
	\end{proposition}
	
	\begin{corollary}
		Let $E(H,K)\equiv \xymatrix{e \ar[r]&H\ar[r]^{i}&G \ar[r]^{\beta}&K\ar[r] & e}$ be a central extension, where $K$ is a finite multiplicative Lie algebra. If $K$ is perfect and $\tilde{M}(K)=e$. Then every pair $(\alpha,\eta) \in Aut(H) \times Aut(K)$ is inducible.
		\begin{proof} \justifying
			Since $K$ is perfect. Let $ E(R,K)\equiv \xymatrix{e \ar[r]&R \ar[r]^{i}&F \ar[r]^{\nu}&K\ar[r] & e}$ be a free presentation of $K$. Then by Proposition \ref{P_1} we have the following universal central extension 
			\[ 	e \parbox{0.8cm}{\rightarrowfill} \frac{([F,F](F \star F))\cap R}{[R,F](R \star F)}  \xrightarrow{\makebox[0.8cm]{$\bar{i}$}} \frac{([F,F](F \star F))}{[R,F](R \star F)} \xrightarrow{\makebox[0.8cm]{$\bar{\nu}$}}   K \longrightarrow  e. \]
			Also, we have $\tilde{M}(K) =e $, this implies 
			\[ 	e \parbox{0.8cm}{\rightarrowfill} e \xrightarrow{\makebox[0.8cm]{$$}} \frac{([F,F](F \star F))}{[R,F](R \star F)} \xrightarrow{\makebox[0.8cm]{$\bar{\nu}$}}   K \longrightarrow  e \]
			is a universal central extension, it can easily seen that $ \xymatrix{e \ar[r]&e \ar[r]^{}&K \ar[r]^{i}& K\ar[r] & e}$ is also universal central extension. Thus by Proposition \ref{P_2}, every central extension by $K$ splits. Hence $H_{ML(\Gamma)}^2(K,H) =e$. Therefore by Theorem \ref{second_second_thm}, every pair $(\alpha,\eta) \in Aut(H) \times Aut(K)$ is inducible.
		\end{proof}
	\end{corollary}
	
	\begin{corollary}
		Let $E(H,K)\equiv \xymatrix{e \ar[r]&H\ar[r]^{i}&G \ar[r]^{\beta}&K\ar[r] & e}$ be a central extension, where $K$ is a cyclic group of order $m$ and $H$ is abelian group of order $n$ such that $m$ and $n$ are coprime. Then every pair $(\alpha,\eta) \in Aut(H) \times Aut(K)$ is inducible.
		\begin{proof}
			We know that the Schur multiplier $M(K)$ of a cyclic group is trivial. Thus by corollary [3.5, \cite{RLS}] we have $H^2_{ML(\Gamma)}(K,H) \cong H^2(K,H) \times Hom(\wedge^2 K,H)$. Since $m$ and $n$ are coprime, $H^2(K,H)=e$. As $M(K)=e$, we have $\wedge^2K=[K,K]=e$. Hence $H^2_{ML(\Gamma)}(K,H) = e$. By Theorem \ref{second_thm}, the map $\Pi$ is surjective.
		\end{proof}
	\end{corollary}
	
	\section{Split Exact Sequences}
	Let $E(H,K)\equiv \xymatrix{e \ar[r]&H\ar[r]^{i}&G \ar[r]^{\beta}&K\ar[r] & e}$ be a center extension. If $C_1^{L^*} =\{\alpha \in C_1^L: \chi_1(\alpha) = e\}$ and $C_2^{L^*} =\{\alpha \in C_2^L: \chi_2(\eta) = e\}$. Then by Theorem \ref{first_thm}, we have the following exact sequences
	\begin{eqnarray}
		\label{eq:1}
		e \parbox{0.8cm}{\rightarrowfill} Aut^{H,K}(G)  \xrightarrow{\makebox[0.8cm]{$i$}} Aut_H^K(G) \xrightarrow{\makebox[0.8cm]{$\Pi_1$}}   C_1^{L^*} \longrightarrow  e \\
		\label{eq:2}
		e \parbox{0.8cm}{\rightarrowfill} Aut^{H,K}(G)  \xrightarrow{\makebox[0.8cm]{$i$}} Aut^H(G) \xrightarrow{\makebox[0.8cm]{$\Pi_2$}}  C_2^{L^*} \longrightarrow e
	\end{eqnarray}
	
	Let $E(H,K)\equiv \xymatrix{e \ar[r]&H\ar[r]^{i}&G \ar[r]^{\beta}&K\ar[r] & e}$ be a central extension, i.e. $H \subseteq Z(G) \cap LZ(G)$ and $C^{L^*} = \{(\alpha,\eta) \in Aut(H) \times Aut(K) : \chi(\alpha,\eta) =e\}$. Then by Theorem \ref{second_thm}, we have the following exact sequence
	\begin{equation}
		\label{eq:3}
		e \parbox{0.8cm}{\rightarrowfill} Z_{ML(\Gamma)}^1(K,H)  \xrightarrow{\makebox[0.8cm]{$i$}} Aut_H(G) \xrightarrow{\makebox[0.8cm]{$\Pi$}}   C^{L^*}  \longrightarrow e
	\end{equation}
	\begin{theorem}
		If exact sequence $E(H,K) \equiv  \xymatrix{	e \ar[r] & H\ar[r]^{i} & G \ar[r]^{\beta} & K\ar[r] & e}$ splits, then the sequences \ref{eq:1} and \ref{eq:2} also split. Also, if $E(H,K)$ is a central extension, i.e. $H\subseteq Z(G) \cap LZ(G)$, then the sequence \ref{eq:3} splits.
		
		\begin{proof}
			Since $E(H,K)\equiv \xymatrix{e \ar[r]&H\ar[r]^{i}&G \ar[r]^{\beta}&K\ar[r] & e}$ splits, we have a multiplicative Lie algebra homomorphism $t: K\to G$ such that $\beta \circ t =I_K$. To show that sequence \ref{eq:1} splits. Define $\delta_1 : C_1^{L^*} \to Aut_H^K(G)$ by $\delta_1(\alpha) = \phi_1$, where $\phi_1: G \to G$ given by $\phi_1(t(x)h) = t(x) \alpha(h)$, $h \in H$ and $x \in K$. Let $g_1=t(x_1)h_1,g_2=t(x_2)h_2 \in G$. Then
			\begin{align*}
				\phi_1(g_1g_2) & =  \phi_1(t(x_1) h_1 t(x_2) h_2)\\
				& = \phi_1(t(x_1x_2) h_1h_2)\\
				& = t(x_1x_2) \alpha(h_1h_2)\\
				& = t(x_1) \alpha(h_1) t(x_2) \alpha(h_2)\\
				& = \phi_1(g_1) \phi_1(g_2).
			\end{align*}
			For the operation $\star$ we have
			\begin{align*}
				\phi_1(g_1 \star g_2) & = \phi_1(t(x_1)h_1 \star t(x_2) h_2)\\
				& = \phi_1((t(x_1) \star t(x_2))( h_1\star t(x_2))(t(x_1)\star h_2))~ (\because ~ H \subseteq Z(G))\\
				& = t(x_1 \star x_2) \alpha(h_1 \star t(x_2)) \alpha(t(x_1)\star h_2) ~ (\because ~ t~ \text{- homomorphism})\\
				& = t(x_1) \star t(x_2) (\alpha(h_1) \star t(x_2)) (t(x_1)\star \alpha(h_2))~ (\because ~ \alpha \in C_1^L) \\
				& = t(x_1) \alpha(h_1) \star t(x_2) \alpha(h_2) \\
				& = \phi_1(g_1) \star \phi_1(g_2).
			\end{align*}
			Thus $\chi_1$ is a multiplicative Lie algebra homomorphism. Let $h \in H$, then clearly $\phi_1(h) =\alpha(h)$, this implies $\phi_1(H)=H$. Consider $\phi_1(g)g^{-1}=\phi_1(t(x)h)(t(x)h)^{-1} = t(x) \alpha(h) (ht(x))^{-1} = \alpha(h) h^{-1} \in H$. Thus $\phi_1 \in Aut_H^K(G)$ and $\Pi_1 \circ \delta_1=I_{C_1^{L^*}}$. \\
			Now we show that the sequence  \ref{eq:2}  splits. Define a map $\delta_2: C_2^{L^*} \to Aut^H(G)$ by $\delta_2(\eta) = \phi_2$, where $\phi_2: G \to G$ given by $\phi(g)= \phi_2(t(x)h) = t(\eta(x)) h$. It can be easily seen that $\phi_2$ is a group homomorphism. Now consider
			\begin{align*}
				\delta_2(g_1 \star g_2) & = \delta_2(t(x_1)h_1\star t(x_2) h_2)\\
				& = \delta_2((t(x_1)\star t(x_2))(t(x_1)\star h_2)(h_1 \star t(x_2)))\\
				& = t(\eta(x_1)\star \eta(x_2)) (t(x_1)\star h_2)(h_1 \star t(x_2)) \\
				& = t(\eta(x_1)) \star t(\eta(x_2)) (t(\eta(x_1))\star h_2) (h_1 \star t(\eta(x_2)))\\
				& = t(\eta(x_1))h_1 \star t(\eta(x_2)) h_2 = \delta_2(g_1) \star \delta_2(g_2).
			\end{align*}
			Thus $\phi_2$ is a multiplicative Lie algebra homomorphism. Clearly $\phi_2(h) = h$, $\phi_2 \in Aut^H(G)$ and $\Pi_2 \circ \delta_2 = I_{C_2^{L^*}}$. Finally for the sequence \ref{eq:3}, define a map $\delta: C^{L^*}: \to Aut_H(G)$ by $\delta(\alpha,\eta) = \phi$, where $\phi: G \to G$ is given by $\phi(g) = \phi(t(x)h) = t(\eta(x)) \alpha(h)$. Since $H \subseteq Z(G) \cap LZ(G)$, i.e. $h \star g = e$ for all $g \in G$ and $h \in H$. Clearly, $\delta$ is a group homomorphism, also we have $\phi(g_1 \star g_2) = \phi(t(x_1 \star x_2))= t(\eta(x_1\ star x_2))= \phi(g_1) \star \phi(g_2)$. Thus $\delta$ is a multiplicative Lie algebra homomorphism. Since $\phi(h) =\alpha(h)$, $\phi \in Aut_H(G)$ and $\Pi \circ \delta = I_{C^{L^*}}$.
	\end{proof}	\end{theorem}

\end{document}